\documentclass[a4paper,12pt]{article}
\usepackage{amscd}
\usepackage{amsmath,amsfonts,amssymb,amscd}
\usepackage{indentfirst,graphicx,epsfig}
\usepackage{graphicx,psfrag}
\usepackage{ifpdf}
\usepackage{float}

\input{epsf}

\newtheorem{thm}{Theorem}

\newtheorem{cor}{Corollary}

\newtheorem{lem}{Lemma}

\def\qed{\hfill \nopagebreak\rule{5pt}{8pt}}

\begin{document}
\title{\Large\bf Rainbow connection number and independence number of
a graph\footnote{Supported by NSFC No.11071130.}}
\author{\small Jiuying Dong, Xueliang Li\\
\small Center for Combinatorics and LPMC-TJKLC\\
\small Nankai University, Tianjin 300071, China\\
\small jiuyingdong@126.com, lxl@nankai.edu.cn}
\date{}
\maketitle

\begin{abstract}

Let $G$ be an edge-colored connected graph. A path of $G$ is called
rainbow if its every edge is colored by a distinct color. $G$ is
called rainbow connected if there exists a rainbow path between
every two vertices of $G$. The minimum number of colors that are
needed to make $G$ rainbow connected is called the rainbow
connection number of $G$, denoted by $rc(G)$. In this paper, we
investigate the relation between the rainbow connection number and
the independence number of a graph. We show that if $G$ is a
connected graph, then $rc(G)\leq 2\alpha(G)-1$. Two examples $G$ are
given to show that the upper bound $2\alpha(G)-1$ is equal to the
diameter of $G$, and therefore the best possible since the diameter
is a lower bound of $rc(G)$. \\

\end{abstract}

\section{Introduction}

All graphs considered in this paper are simple, finite and
undirected. The following notation and terminology are needed in the
sequel. Let $u\in V$ and $v\in V$ be two distinct vertices of a
graph $G=(V,E)$. The \emph{distance between $u$ and $v$} in $G$,
denoted by $d(u,v)$, is the length of a shortest path connecting
them in $G$. A \emph{$(u,v)$-path} is a path with initial vertex $u$
and terminal vertex $v$, denoted by $P[u,v]$. Let $P_H[u,v]$ denote
the path in $H$ connecting $u$ and $v$, where $H$ is the subgraph of
$G$. For two subsets $U$ and $W$ of $V$, a $(U,W)$-path is a path
which starts at a vertex of $U$ and ends at a vertex of $W$, and
whose internal vertices belong to neither $U$ nor $W$. We use
$E[U,W]$ to denote the set of edges of $G$ with one end in $U$ and
the other end in $W$, and $e(U,W)=|E[U,W]|$. Let $G[U]$ denote the
subgraph of $G$ whose vertex set is $U$ and whose edge set consists
of all such edges of $G$ that have both ends in $U$. The following
notions are from \cite{K-Y}. A set $D\subseteq V(G)$ is called a
$k$-step connected dominating set of $G$, if every vertex in
$G\setminus D$ is at a distance at most $k$ from $D$, where $G[D]$
is connected. The $k$-step open neighborhood of a set $D$ is
$N^{k}(D):=\{x\in V(G)|d(x,D)=k\}$, where $k\in N$. We often use
$e(G)$ to denote the number of edges in a graph $G$ and $|G|$ to
denote the order of $G$. For undefined terminology and notation, we
refer to \cite{B-M}.

Let $G=(V,E)$ be a connected graph with vertex set $V$ and edge set
$E$.  A \emph{$k$-edge coloring }of $G$ is a mapping $c: E
\rightarrow C$, where $C$ is a set of $k$ distinct colors. In
\cite{C-J} Chartrand, Johns, McKeon and Zhang introduced a new
concept about the connectivity and coloring of a graph, which is
given follows. A path of $G$ is called \emph{rainbow} if every edge
of it is colored by a distinct color. For every two vertices $u$ and
$v$ of $G$, if there exists a rainbow path between them, we say that
$G$ is \emph{rainbow connected}. \emph{The rainbow connection number
$rc(G)$} is defined as the smallest number of colors that are needed
to make $G$ rainbow connected. An edge coloring is called \emph{a
rainbow coloring} if it makes $G$ rainbow connected. From the
definition of rainbow connection, we can see that the diameter
$diam(G)\leq rc(G)\leq e(G)$. For more knowledge on the rainbow
connection, we refer to \cite{L-S1, L-S2}.

In \cite{Chakraborty} Chakraborty, Fischer, Matsliah and Yuster
showed that given a graph $G$, deciding if $rc(G)=2$ is NP-complete,
in particular, computing $rc(G)$ is NP-hard, which were conjectured
by Caro, Lev, Roditty, Tuza and Yuster \cite{Caro}. So, to obtain
upper bounds for the rainbow connection number $rc(G)$ of a graph
$G$ becomes interesting. Therefore, many good upper bounds have been
obtained in terms of other graph parameters. Caro, Lev, Roditty,
Tuza and Yuster\cite{Caro} conjectured that if $G$ is a connected
graph with $n$ vertices and $\delta(G)\geq 3$, then $rc(G)<
\frac{3}{4}n$. Schhiermeyer \cite{Schiermeyer2} confirmed the
conjecture and showed that if $G$ is a connected graph with $n$
vertices and $\delta(G)\geq 3$, then $rc(G)\leq\frac{3n-1}{4}$, and
there are examples to show that $\frac{3}{4}$ cannot be replaced
with a smaller constant. In \cite{Chandran} Chandran, Das,
Rajendraprasad and Varma obtained a good relation between the
rainbow connection number and the minimum degree of a graph. They
showed that if $G$ is a connected graph of order $n$ and minimum
degree $\delta(G)$, then $rc(G)\leq 3n/(\delta(G)+1)+3$, and the
bound is tight up to addictive factors. Later, we \cite{Dongli}
studied the relation between the rainbow connection number and the
minimum degree sum, a generalized result of the minimum degree
version. We showed that if $G$ is a graph with $k$ independent
vertices, then $rc(G)\leq\frac{3kn}{\sigma_{k}(G)+k}+6k-3$. In
\cite{Basavaraju}, Basavaraju, Chandran, Rajendraprasad and
Ramaswamy investigated the relation between the rainbow connection
number and the radius of a bridgeless graph. They showed that for
every bridgeless graph $G$ with radius $rad(G)$, $rc(G)\leq
rad(G)(rad(G)+2)$, and gave an example to show that the bound is
tight. In \cite{Li} Li, Liu, Chandran, Mathew and Rajendraprasad
showed that if $G$ is a 2-connected graph of order $n \ (n\geq 3)$,
then $rc(G)\leq\lceil\frac{n}{2}\rceil$, and the upper bound is
tight for $n\geq 4$. Later, Ekstein et al. \cite{Eks} got the same
result. Li et al. \cite{Li} also got some relations between the
rainbow connection number and the connectivity of a graph.
Schiermeyer \cite{Sch1} obtained a relation between the rainbow
connection number of a graph $G$ and the chromatic number of the
complement of $G$, i.e., $rc(G)\leq 2\chi(\bar{G})-1$.

This paper intends to give a relation between the rainbow connection
number and the independence number of a graph. Recall that an
\emph{independent set} of a graph $G$ is a set of vertices such that
any two of these vertices are non-adjacent in $G$. \emph{The
independence number} $\alpha(G)$ of $G$ is the cardinality of a
maximum independent set of $G$. The independence number of a graph
is an important parameter, and the investigation on the relations
between the independence number and other graph parameters is
interesting, see \cite{Brandt, Chen, Harant, Lichiardopol, Luo3,
Luo1, Luo2}. We obtain the following result.
\begin{thm}
If $G$ is a connected graph, then $rc(G)\leq 2\alpha(G)-1$,
and the bound is the best possible.
\end{thm}

Two examples are given to show that our bound $2\alpha(G)-1$ is
exactly equal to the diameter of $G$, and therefore our bound is the
best possible since the diameter is a lower bound of $rc(G)$.
Moreover, for these examples some good bounds in terms of other
parameters can be arbitrarily bad. As we know, even for the
chromatic number, every upper bound is attained yet arbitrarily bad
for many graphs. We also give an example, Example 3, to show that
our result $rc(G)\leq 2\alpha(G)-1$ could be arbitrarily bad. This
example also shows that the bounds in terms of other parameters are
arbitrarily bad.

{\noindent\bf Example 1}: Let $P_{2t}=v_1v_2v_3\cdots
v_{2t-1}v_{2t}$ be a path of length $2t-1$, and let $G_1, G_2,
\cdots, G_t$ be $t \ (t\geq 2)$ complete graphs with $|G_1|=2$ and
$|G_i|=s$ for $2\leq i\leq t$. For every $i$ with $1\leq i\leq t$,
we join each vertex of $G_i$ to every vertex of $v_{2i-1}$ and
$v_{2i}$. The obtained graph is denoted by $G$. One can see that $G$
is connected with $\delta(G)=3$, and
$I(G)=\{v_2,v_4,v_6,\cdots,v_{2t}\}$ is a maximum independent set,
that is, $\alpha(G)=t$. We also know that the distance
$d(v_1,v_{2t})=2t-1$. So, $rc(G)\geq 2t-1$. Now we use $2t-1$
distinct colors to give $G$ an edge coloring. Let $1,2,\cdots,2t-1$
be $2t-1$ distinct colors. We use the $2t-1$ colors to color all the
edges of $P_{2t}$, and each edge with a distinct color. Thus,
$P_{2t}$ is rainbow connected. Then we use color $2i-1$ to color
every edge of $E[V(G_i),\{v_{2i-1},v_{2i}\}]$. Finally, we use color
1 to color every edge of $G[V(G_i)]$. One can show that $G$ is
rainbow connected. For each pair $(u,v)\in(G_i, P_{2t})$, either the
edge $uv_{2i}$ together with the path in $P_{2t}$ connecting
$v_{2i}$ and $v$ forms a rainbow path, or the edge $uv_{2i-1}$
together with the path in $P_{2t}$ connecting $v_{2i-1}$ and $v$
forms a rainbow path. For each pair $(u,v)\in(G_i, G_j)$ with $1\leq
i < j\leq t$, the edges $uv_{2i}$ and $vv_{2j-1}$ together with the
path in $P_{2t}$ connecting $v_{2i}$ and $v_{2j-1}$ form a rainbow
path. So, $G$ is rainbow connected, and hence $rc(G)\leq 2t-1$.
Thus, $rc(G)=2t-1=2\alpha(G)-1$. Note that, $diam(G)=2t-1$,
$rad(G)=t$, $\delta(G)=3$ and for any $v\in
V(G)\setminus(\{v_1,v_2\}\cup V(G_1))$, the degree of $v$ is at
least $s+1$.

The following facts can be easily seen. When $s$ is very large, the
order $n$ of $G$ is very large. The upper bounds in \cite{Chandran,
Schiermeyer2} can be arbitrarily bad, because $
3n/(\delta(G)+1)+3=3n/4+3$ and $\frac{3n-1}{4}$, both bounds are
large by selecting $s$ to be large; When $k\geq 2$, $\sigma_k(G)\geq
3+s+1$, the bound in \cite{Dongli} is
$\frac{3kn}{\sigma_k(G)+k}+6k-3\leq\frac{3kn}{s+4+k}+6k-3$, better
than the bounds of \cite{Chandran, Schiermeyer2}, when $s$ is very
large, but it is also far from the diameter of $G$; The bound in
\cite{Basavaraju} is $rad(G)(rad(G)+2)=t(t+2)$, a square of $t$.
However, our result $rc(G)\leq 2\alpha(G)-1$ is the best, which is
equal to $diam(G)$.

{\noindent\bf Example 2}: For $1 \leq i\leq 2t \ (t\geq 2)$, let the
sets $V_1, V_2, \cdots, V_{2t}$ be pairwise disjoint,
$|V_1|=|V_2|=2$ and  for $i\geq 3$, $|V_i|=s \ (s\geq 2)$. Join each
vertex of $V_i$ to every vertex of $V_{i+1}$ for $1\leq i\leq 2t-1$.
Suppose that every two vertices of $V_i$ are adjacent for each $i$.
We denote the resulting graph by $G$. Note that, $G$ is 2-connected,
$\alpha(G)=t$, $diam(G)=2t-1$ and $rad(G)=t$. For any vertex $v\in
V_1$, the degree of $v$ is 3, for any vertex $v\in V_2$, the degree
of $v$ is $s+3$, for any vertex $v\in V\setminus\{V_1,V_2,
V_{2t}\}$, the degree of $v$ is $3s-1$, and for any vertex $v\in
V_{2t}$, the degree of $v$ is $2s-1$. Color each edge in $E[V_i,
V_{i+1}]$ with color $i$ for $1\leq i\leq 2t-1$ and each edge in
$G[V_i]$ with color 1 for $1 \leq i\leq 2t$. It is not difficult to
check that $rc(G)=2\alpha(G)-1=2t-1=diam(G)$.

One can see the following facts.  When $s$ is very large, the order
$n$ of $G$ is very large. The upper bounds in \cite{Chandran, Li,
Schiermeyer2} can be arbitrarily bad, because $
3n/(\delta(G)+1)+3=3n/4+3$, $\lceil\frac{n}{2}\rceil$ and
$\frac{3n-1}{4}$ are all large by selecting $s$ to be large; When
$s$ is very large, the bound in \cite{Dongli} is
$\frac{3kn}{\sigma_k(G)+k}+6k-3\leq\frac{3kn}{s+6+k}+6k-3$, better
than the bounds of \cite{Chandran, Li Schiermeyer2}, but it is also
far from $diam(G)$. The bound in \cite{Basavaraju} is
$rad(G)(rad(G)+2)=t(t+2)$, a square of $t$, far from $diam(G)$.
However, our result $rc(G)\leq 2\alpha(G)-1$ is the best, which is
equal to $diam(G)$.

{\noindent\bf Example 3}: Consider the graph $G=K_{1,1,s}$ with
partition sets $V_1$, $V_2$ and $V_3$, and $|V_1|=|V_2|=1$, and
$|V_3|=s \ (s > 3)$. Then $G$ is 2-connected, $\delta(G)=2$,
$\alpha(G)=s$. and $rc(G)=3$. We can see that the upper bounds of
\cite{Chandran, Li, Schiermeyer2} and our result can be arbitrarily
bad, because $3n/(\delta(G)+1)+3=3n/4+3$, $\frac{3n-1}{4}$,
$\lceil\frac{n}{2}\rceil$ and $2s-1$ are very large by selecting $s$
to be large.

\section{Proof of Theorem 1}

Before proving our main result, we first prove a lemma. Although
this lemma can be found in \cite{Chandran}, here we use a new
technique to give it another proof.

\begin{lem}
Let $G$ be a connected graph with $\delta(G)\geq 2$, and let $D$ be
a connected dominating set of $G$. Then $rc(G)\leq rc(G[D])+3$.
\end{lem}

\noindent{\bf Proof.} At first, we use $rc(G[D])$ different colors
to give $G[D]$ a rainbow coloring. Then, let $1,2$ and $3$ be three
distinct fresh colors. Since $D$ is a connected dominating set of
$G$, $V(G)=D\cup N(D)$. We will perform the following procedure:

\noindent{\bf Procedure 1:}

 \ \ \ \ \ \ $F= N(D)$, $i=1$,\\
\hspace*{\fill} while there exists a vertex $w_i\in F$ with $d_{G[F]}(w_i)\geq 1$
do \hspace*{\fill}\\
 \hspace*{\fill}   $H_i=N[w_i]\cap F$, $F=F\setminus H_i$,\hspace*{\fill}\\
 \hspace*{\fill}   $i=i+1$.\hspace*{\fill}\\
While the above procedure ends, we get an empty graph $G[F]$, and a
sequence of vertices $w_1,w_2,\cdots, w_t$ and a sequence of sets
$H_1, H_2,\cdots, H_t$. So, we have partitioned $N(D)$ into some
disjoint subsets $H_1, H_2,\cdots, H_t, F$. Now we will give colors
to the remaining uncolored edges of $G$. We use color $1$ to color
every edge in $E[w_i,D]$, and use color $2$ to color every edge in
$G[N(D)]$. Finally, we use color 3 to color every edge in
$E[H_i\setminus\{w_i\}, D]$. It is not difficult to check that $G$
is rainbow connected. \qed

Before giving the proof of Theorem 1, we need the following
observation.

\noindent{\bf Observation.} Let $G$ be a graph. If $G$ has a cut
edge $uv$, then we replace $uv$ by a clique of order at least 3,
i.e., we add to $G$ some new vertices $w_1$, $w_2\ldots, w_q$ with
$q\geq 1$ such that $u,v, w_1, \ldots, w_q$ form a complete
subgraph. The new graph is denoted by $G'$, and is called a {\it
blow-up graph} of $G$ at the cut edge $uv$. It is not difficult to
check that $rc(G')=rc(G)$ and $\alpha(G')=\alpha(G)$. In this way,
we need only to consider graphs without any cut edge, and therefore
without any pendant edge, i.e., its minimum degree is at least 2.

{\noindent\bf Proof of Theorem 1.} If $G$ is a complete graph, then
$\alpha(G)=1$ and $rc(G)=1$, Theorem 1 follows. Now assume that $G$
is a non-complete graph with $\delta(G)\geq 2$.

We will perform the following procedure to obtain
a tree $T$ whose vertex set $D$ is a connected dominating set of $G$.
Let $y_0\in V(G)$ with $d(y_0)=\delta(G)$. Since $G$ is a
non-complete graph, $N^{2}(y_0)\neq\emptyset$. We look at the
following procedure:

\noindent {\bf Procedure 2:}

 \ \  $D=\{y_0\}, T=y_0, X=\phi, Y=\{y_0\}.$\\
\ \ \ \ \  While $N^{2}(D)\neq\phi$\\
\hspace*{\fill} take any vertex $v\in N^{2}(D)$, let $P=vhu $ be a path of length 2, \hspace*{\fill}\\
\hspace*{\fill}  where $h\in N^{1}(D)$ and $u\in D$. Let $D=D\cup V(P)$,\hspace*{\fill}\\
\hspace*{\fill}   $ T=T\cup P$, $X=X\cup\{h\}, Y=Y\cup\{v\}$.\hspace*{\fill}\\
If $u\in X$, we call $u$ an $X$-knot vertex. Note that $N^{2}(D)$
does not contain any neighbor of $Y$. When the above procedure ends,
the algorithm runs $|X|$ rounds. Thus we get $V(G)=D\cup N^{1}(D)$,
where $D$ is a connected dominating set. Note that $Y$ is an
independent set and $|Y|=|X|+1$. So $|Y|\leq\alpha(G)$ and
$|D|=|Y|+|X|=2|Y|-1$. Note that $T$ is a spanning tree of $G[D]$ and
the pendant vertices of $T$ are all in $Y$, and also note that if
$x$ is an $X$-knot vertex, then $x$ is adjacent to at least three
vertices of $T$.

In the following, $D$, $T$, $Y$ and $X$ are always the same as those
obtained in the above algorithm. In order to continue our proof, we
need the following lemmas.
\begin{lem}
If there exists a vertex $w\in N^{1}(D)$ such that $e(w,Y)=0$, then
$rc(G)\leq 2\alpha(G)-1$.
\end{lem}
\noindent{\bf Proof}. Let $I=Y\cup\{w\}$. Then $I$ is an independent
set and $|I|=|Y|+1$. So $|Y|=|I|-1\leq\alpha(G)-1$. By Lemma 1, we
can get that $rc(G)\leq rc(G[D])+3\leq |D|+2=2|Y|+1$. Hence,
$rc(G)\leq 2(\alpha(G)-1)+1=2\alpha(G)-1$. \qed

From the proof of Lemmas 2, we can see that if we can find an independent
set $I$ satisfying $|I|=|Y|+1$, then we will get $rc(G)\leq 2\alpha(G)-1$.

\begin{lem}
If $G[D]=T$ and there exist two vertices $w, w'\in N^{1}(D)$ such
that $ww'\not\in E(G)$, $e(w,Y)=1$, $e(w',Y)=1$, $e(w,X)=0$ and
$e(w',X)=0$, then $rc(G)\leq 2\alpha(G)-1$.
\end{lem}

\noindent{\bf Proof.} Let $wy\in E(G)$ and $w'y'\in E(G)$ where
$y\in D$ and $y'\in D$. Since $G[D]=T$, let $P_T[y,y']$ denote the
only path connecting $y$ and $y'$ in $T$. If there do not exist two
successive vertices of $X$ in $P_T[y,y']$, then let
$I=\{w,w'\}\cup(V(P_T[y,y'])\cap X)\cup(Y\setminus(V(P_T[y,y'])\cap
Y))$. One can see that $I$ is an independent set and $|I|=|Y|+1$,
and so $rc(G)\leq 2\alpha(G)-1$. Otherwise, suppose that there exist
two successive vertices of $X$ in $P_T[y,y']$. By the structure of
$T$, we can conclude that there must be an $X$-knot vertex between
the two successive vertices. Then there must be a segment in
$P_T[y,y']$, without loss of generality, say $P_T[y,x]\subset
P_T[y,y']$, where $x$ is an $X$-knot vertex, and in $P_T[y,x]$ there
is an vertex $x'$ of $X$ adjacent to $x$, and $x',x$ are the only
two successive vertices in $P_T[y,x]$. Then
$I=\{w\}\cup(V(P_T[y,x'])\cap X)\cup(Y\setminus V(P_T[y,x'])\cap Y)$
is an independent set and $|I|=|Y|+1$. So $rc(G)\leq
2\alpha(G)-1$.\qed

\begin{lem}
If $G[D]=T$ and there exist vertices $w_1, w_2\in N^{1}(D)$ and
$y_1, y_2\in D$ such that $w_1w_2\not\in E(G)$, and $N(w_1)\cap
D=N(w_2)\cap D=\{y, y'\}$, then $rc(G)\leq 2\alpha-1$.
\end{lem}

\noindent{\bf Proof.} Let $P_T[y,y']$ denote the only path
connecting $y$ and $y'$ in $T$. If there do not exist two successive
vertices of $X$ in $P_T[y,y']$, then let $I=\{w_1,
w_2\}\cup(V(P_T[y,y'])\cap X)\cup(Y\setminus(V(P_T[y,y'])\cap Y))$.
One can see that $I$ is an independent set and $|I|=|Y|+1$, and so
$rc(G)\leq 2\alpha(G)-1$. Otherwise, there must exist two successive
vertices of $X$ in $P_T[y,y']$. By the structure of $T$, we can
conclude that there must be an $X$-knot vertex between the two
successive vertices. Then there must be a segment in $P_T[y,y']$,
without loss of generality, say $P_T[y,x]\subset P_T[y,y']$, where
$x$ is an $X$-knot vertex, and in $P_T[y,x]$ there is an vertex $x'$
of $X$ adjacent to $x$, and $x',x$ are the only two successive
vertices in $P_T[y,x]$. Then $I=\{w,w'\}\cup(V(P_T[y,x'])\cap
X)\cup((Y\setminus((V(P_T[y,x']) \cap Y)\cup\{y'\})$ is an
independent set and $|I|=|Y|+1$. So $rc(G)\leq 2\alpha(G)-1$.\qed

Let $N^{1}(D)=A\cup B$ where $w\in A$ if and only if $e(w,D)\geq 2$,
and $w\in B$ if and only if $e(w,D)=1$. By Lemma 2, we can get that
every vertex $w\in B$ satisfies $e(w,X)=0$, and so $e(w,Y)=1$. By
Lemma 3, we can get that $G[B]$ is a complete subgraph. In the
following we will divide two cases to finish our proof.

\noindent{\bf Case 1.} $e(G[D])\geq e(T)+1$.

Let $a_1a_2\in E(G[D])$ and $a_1a_2\not\in E(T)$. Note that $T$ is a
spanning tree of $G[D]$. So $T\cup a_1a_2$ contains a cycle, say
$C$, and $a_1a_2\in E(C)$. Let $G'=T\cup a_1a_2$. Then $rc(G[D])\leq
rc(G')$. Since $rc(G')\leq e(T)-(|C|-1)+rc(C)$ and $rc(C)\leq
\lceil\frac{|C|}{2}\rceil$ when $|C|\geq 4$, we can get

\[
  rc(G')\leq\left\{
  \begin{array}{ll}
  e(T)-\frac{|C|}{2}+1, & |C| \  is\  even \\
  e(T)-\frac{|C|-3}{2}, & |C|\  is \  odd \  and \ |C|\neq 3 \\
  e(T)-1, & |C|= 3
  \end{array}
  \right.
\]
Hence, $rc(G[D])\leq rc(G')\leq e(T)-1$.

Now we color every edge of $G$ in the following way. First, we use
$rc(G[D])$ distinct colors to rainbow color $G[D]$. Then let $c',
c''$ be two fresh colors. For any vertex $w\in A$, let $w', w'' \in
D$ with $ww', ww''\in E(G)$, set $c(ww')=c'$ and $c(ww'')=c''$.
For any vertex $b\in B$, let $b'\in D$ with $bb'\in E(G)$, set
$c(bb')= c'$. For the remaining uncolored edges of $E(G)$, we use a
used color to color them. Thus we have colored all the edges of $G$.

We will show that $G$ is rainbow connected. For each pair
$(u,v)\in(N(D)\times D)$, the edge $uu'$ together with the path in
$G'$ connecting $u'$ and $v$ forms a rainbow path, where
$c(uu')=c'$ and $u'\in D$. For each pair $(u,v)\in(A\times A)$,
the edges $uu'$ and $vv''$ together with the path in $G'$
connecting $u'$ and $v''$ form a rainbow path, where $c(uu')=c'$
and $c(vv'')=c''$. For each pair $(u,v)\in(A\times B)$, the edges
$uu''$ and $vv'$ together with the path in $G'$ connecting $u''$
and $v'$ form a rainbow path, where $c(uu'')=c''$ and $c(vv')=c'$.
Thus we have showed that $G$ is rainbow connected.

In the above coloring, we used at most $rc(G[D])+2\leq e(T)+1$
colors. Hence, $rc(G)\leq e(T)+1$, that is $rc(G)\leq |D|$. Since
$|D|=2|Y|-1\leq 2\alpha(G)-1$, we can get $rc(G)\leq 2\alpha(G)-1$.

\noindent {\bf Case 2.} $e(G[D])=e(T)$.

Let $1,2,c_1,c_2$ and $a$ be 5 distinct colors, and in the following
proof, we use $a$ to color each edge of $E[B,D]$,  and use $c_1$ to color each
edge of $ E(G[B])$.

Choose a longest path $P$ in $G[D]$ such that two ends of $P$ are
pendant vertices. We know that the two pendant vertices belong to
$Y$, and $|P|\geq 3$. Let $P=y_1x_1x_2\cdots x_ky_2$. We look at two
subcases:

\noindent {\bf Subcase 2.1.} $V(P)\subset D$.

Since $P$ is a longest path and $|Y|=|X|+1$, we can get $|P|\geq 4$.
Choose a pendant edge in $T$, say $y_3x$, which is not in $P$. Let
$P'$ be a path passing through $y_3x$ in $T$ with $|V(P)\cap V(P')|=1$,
and let $V(P)\cap V(P')=\{x'\}$. Without loss
of generality, let $|P[y_1,x']|\geq 3$, and let $c(y_1x_1)=1$,
$c(x_1x_2)=c_1$, $c(x_ky_2)=2$ and $c(xy_3)=c_2$.

Let $A_1$, $A_2$, $A_3$ and $A_4$ be the subsets of $A$. $w_1\in
A_1$ if and only if $w_1y_1\in E(G)$ and $w_1$ is adjacent to only
one vertex of $D\setminus\{y_1, y_2\}$. Let $c(w_1y_1)=c_2$ and
$c(w_1w_1')=2$ where $w_1'\in D$; $w_2\in A_2$ if and only if
$w_2y_2\in E(G)$ and $w_2$ is adjacent to only one vertex of
$D\setminus\{y_1, y_2\}$. Let $c(w_2y_2)=c_2$ and $c(w_2w_2')=1$
where $w_2'\in D$; $w_3\in A_3$ if and only if $w_3y_1\in E(G)$,
$w_3y_2\in E(G)$ and $e(w_3, D)=2$. Let $c(w_3y_1)=2$ and
$c(w_3y_2)=1$. $w_4\in A_4$ if and only if $w_4$ is adjacent to at
least two vertices $w_4'$ and $w_4''$ of $D\setminus\{y_1, y_2\}$.
Assume that the distance between $w_4'$ and $y_1$ in $T$ is not
more than the distance between $w_4''$ and $y_1$ in $T$. Let
$c(w_4w_4')=a$ and $c(w_4w_4'')=1$. Let $B_1$, $B_2$ and $B_3$ be
the subsets of $B$. $b_1\in B_1$ if and only if $b_1$ is only
adjacent to $y_1$; $b_2\in B_2$ if and only if $b_2$ is only
adjacent to $y_2$; $b_3\in B_3$ if and only if $b_3$ is only
adjacent to some vertex of $Y\setminus\{y_1,y_2\}$.
Thus we get $A=A_1\cup A_2\cup A_3\cup A_4$, $B=B_1\cup B_2\cup B_3$.
From Lemma 4 we know that $G[A_3]$ is a
complete subgraph, and from Lemma 3 we get that $G[B]$ is a complete
subgraph. It is easy to check that for any vertex of $N(D)$ is
rainbow connection to any vertex of $D$.

\noindent {\bf Subsubcase 2.1.1.} $B_1\neq\phi$, $B_2\neq\phi$ and $B_3\neq\phi$.

First, we will show that $G[A_1]$, $G[A_2]$ and $G[A_4]$ is rainbow
connected, respectively. For each pair $(u,v)\in(A_1\times A_1)$,
let $u',v'\in D$ with $uu',vv'\in E(G)$, if $u'\neq v'$, then without loss of generality,
we assume that the path in $T$ from $v'$ to $y_1$ does not contain the edge $y_3x$.
Thus the edges $uy_1$ and $vv'$ together with the path in $T$
connecting $y_1$ and $v'$ form a rainbow path between $u$ and $v$;
if $u'=v'$ and $v'y_2 \in E(G)$, then the edges $uy_1$ and $vy_2$
together with the path $P$ form a rainbow path between $u$ and $v$;
if $u'=v'$ and $u'y_2\in E(G)$,
similarly, there is a rainbow path between them; if $u'=v'$ and
assume that $v'y_2 \not\in E(G)$ and $u'y_2 \not\in E(G)$, then
from Lemma 4, we can get $uv\in E(G)$. So for each pair $(u,v)\in(A_1\times A_1)$,
there is a rainbow path connecting them.
For each pair $(u,v)\in(A_2\times A_2)$, similar to $(u,v)\in(A_1\times A_1)$, we
can get a rainbow path between $u$ and $v$; For each pair
$(u,v)\in(A_4\times A_4)$, the edges $uu'$ and $vv'$ together with
the path in $T$ connecting $u'$ and $v'$ form a rainbow path,
where $c(uu')=a$ and $c(vv')=1$.

Second, we will show that for any vertex $u\in A_1$, there is a
rainbow path connecting it to any vertex of $A_2\cup A_3\cup A_4\cup
B$. For each pair $(u,v)\in(A_1\times(A_2\cup A_4))$, the edges
$uu'$ and $vv'$ together with the path in $T$ connecting $u'$ and
$v'$ form a rainbow path, where $c(uu')=2$ and $c(vv')=1$. For
each pair $(u,v)\in(A_1\times A_3)$, $uy_1v$ is a rainbow path. For
each pair $(u,v)\in(A_1\times B_1)$, $uy_1v$ is a rainbow path. For
each pair $(u,v)\in(A_1\times B_2)$, the edges $uy_1$ and $vy_2$
together with the path $P$ form a rainbow path. For each pair
$(u,v)\in(A_1\times B_3)$, the edges $uu'$ and $vv'$ together with
the path in $T$ connecting $u'$ and $v'$ form a rainbow path,
where $c(uu')=2$ and $c(vv')=a$.

Third, we will show that for any vertex $u\in A_2$, there is a
rainbow path connecting it to any vertex of $A_3\cup A_4\cup B$. For
each pair $(u,v)\in(A_2\times A_4)$, the edges $uu'$ and $vv'$
together with the path in $T$ connecting $u'$ and $v'$ form a
rainbow path, where $c(uu')=1$ and $c(vv')=a$. For each pair
$(u,v)\in(A_2\times A_3)$, $uy_2v$ is a rainbow path. For each pair
$(u,v)\in(A_2\times B_1)$, the edges $uy_2$ and $vy_1$ together with
the path $P$ form a rainbow path. For each pair $(u,v)\in(A_2\times
B_2)$, $uy_2v$ is a rainbow path. For each pair $(u,v)\in(A_2\times
B_3)$, the edges $uu'$ and $vv'$ together with the path in $T$
connecting $u'$ and $v'$ form a rainbow path, where $c(uu')=1$
and $c(vv')=a$.

Fourth, we will show that for any vertex $u\in A_3$, there is a
rainbow path connecting it to any vertex of $A_4\cup B$. For each
pair $(u,v)\in(A_3\times A_4)$, the edges $uy_1$ and $vv'$ together
with the path in $T$ connecting $y_1$ and $v'$ form a rainbow path,
where $c(uy_1)=2$ and $c(vv')=a$. For each pair $(u,v)\in(A_3\times
B_1)$, $uy_1v$ is a rainbow path. For each pair $(u,v)\in(A_3\times
B_2)$, $uy_2v$ is a rainbow path. For each pair $(u,v)\in(A_3\times
B_3)$, the edges $uy_1$ and $vv'$ together with the path in $T$
connecting $y_1$ and $v'$ form a rainbow path.

Finally, we will show that for any vertex $u\in A_4$, there is a
rainbow path connecting it to any vertex of $B$. For each pair
$(u,v)\in(A_4\times B_1)$, the edges $uu'$, $vb_2$ and $b_2y_2$
together with the path in $T$ connecting $u'$ and $y_2$ form a
rainbow path, where $c(uu')=1$ and $b_2\in B_2$. For each pair
$(u,v)\in(A_4\times B_2)$, the edges $uu'$ and $vy_2$ together with
the path in $T$ connecting $u'$ and $y_2$ form a rainbow path,
where $c(uu')=1$. For each pair $(u,v)\in(A_4\times B_3)$, the
edges $uu'$ and $vv'$ together with the path in $T$ connecting
$u'$ and $v'$ form a rainbow path, where $c(uu')=1$ and
$c(vv')=a$.

Thus, we have proved that $G$ is rainbow connected.

From the proof above, we can see the following facts: for any vertex
of $A$ there is a rainbow path connecting it to any vertex of $G$,
and the internal vertex of the rainbow path does not contain any
vertex of $B$; for any vertex of $B_2$, there is a rainbow path
connecting it to any vertex of $G$, and the rainbow path does not
contain any vertex of $B_1\cup B_3$; for any vertex of $B_3$, there
is a rainbow path connecting it to any vertex of $G$, and the
rainbow path does not contain any vertex of $B_1\cup B_2$. Hence, in
the following proof, we can assume that $B_3=\phi$ and $B_2=\phi$.

\noindent {\bf Subsubcase 2.1.2.} $B_1\neq\phi$.

We still use the above mentioned subsets $A_1$, $A_2$, $A_3$, $A_4$
and $B_1$, and we still color the edges of $G$ in the above way
except for setting $c(w_4w_4')=a$ and $c(w_4w_4'')=2$. Thus we
only need to show that for any vertex of $A_4$, there is a rainbow
path connecting it to any vertex of $G$. We will give the proof as
follows. For each pair $(u,v)\in(A_4\times A_4)$, the edges $uu'$
and $vv'$ together with the path in $T$ connecting $u'$ and $v'$
form a rainbow path, where $c(uu')=a$ and $c(vv')=2$. For each
pair $(u,v)\in(A_4\times A_3)$, the edges $uu'$ and $vy_1$ together
with the path in $T$ connecting $u'$ and $y_1$ form a rainbow path,
where $c(uu')=a$. For each pair $(u,v)\in(A_4\times A_2)$, the
edges $uu'$ and $vv'$ together with the path in $T$ connecting
$u'$ and $v'$ form a rainbow path, where $c(uu')=a$ and
$c(vv')=1$. For each pair $(u,v)\in(A_4\times A_1)$, the edges
$uu'$ and $vv'$ together with the path in $T$ connecting $u'$ and
$v'$ form a rainbow path, where $c(uu')=a$ and $c(vv')=2$. For
each pair $(u,v)\in(A_4\times B_1)$, the edges $uu'$ and $vy_1$
together with the path in $T$ connecting $u'$ and $y_1$ form a
rainbow path, where $c(uu')=2$. Hence, we have showed that $G$ is
rainbow connected.

\noindent {\bf Subcase 2.2.} $V(P)=D$.

Since $V(P)=D$ and $|Y|=|X|+1$, $P$ is a $(Y,X)$-alternate path. Let
$A_1$, $A_2$, $A_3$, $A_4$, $B_1$, $B_2$ and $B_3$ be the above
mentioned subsets.

\noindent {\bf Subsubcase 2.2.1.} $|P|=3$.

Let $P=y_1x_1y_2$. We use color 1 to color edge $y_1x_1$ and use color 2 to color $y_2x_1$.
Let $A_1$, $A_2$, $A_3$, $B_1$ and $B_2$ be the above
mentioned subsets. Note that: $A_4=\phi$, $B_3=\phi$,
$G[A_1\cup B_1]$ is a complete subgraph, and $G[A_2\cup B_2]$ is a complete subgraph.
For any $w_1\in A_1$ and $w_2\in A_2$, let $c(w_1y_1)=a$, $c(w_1x_1)=1$,
$c(w_2y_2)=a$ and $c(w_2x_1)=2$. It is obvious that for each vertex of $A\cup B$, there
is a rainbow path connecting it to any vertex of $P$.
For each pair $(u,v)\in(A_1\times A_2)$,
$ux_1v$ is a rainbow path. For each pair $(u,v)\in(A_1\times A_3)$,
$uy_1v$ is a rainbow path. For each pair $(u,v)\in(A_1\times B_2)$,
$ux_1y_2v$ is a rainbow path.  For each pair $(u,v)\in(A_2\times A_3)$,
$uy_2v$ is a rainbow path. For each pair $(u,v)\in(A_2\times B_1)$,
$ux_1y_1v$ is a rainbow path. For each pair $(u,v)\in(A_3\times B_1)$,
$uy_1v$ is a rainbow path. For each pair $(u,v)\in(A_3\times B_2)$,
$uy_2v$ is a rainbow path. Thus, we have showed that $G$ is rainbow connected.

\noindent {\bf Subsubcase 2.2.2.} $|P|\geq 5$.

Let $c(y_1x_1)=1$, $c(x_1y_1')=c_1$,
$c(y_2x_2)=2$ and $c(x_2y_2')=c_2$ where $y_1', y_2'\in V(P)$. We
color the edges of $G$ in the following way: We use $a$ to color each edge of
$E[B,D]$, and use $c_1$ to color each edge of $G[B]$.
For any $w_1\in A_1$, let
$c(w_1y_1)=2$ and $c(w_1w_1')=a$ where $w_1'\in D$; For any $w_2\in A_2$,
let $c(w_2y_2)=1$ and $c(w_2w_2')=a$ where $w_2'\in D$; For any
$w_3\in A_3$, let $c(w_3y_1)=2$ and $c(w_3y_2)=1$; For any $w_4\in
A_4$, assume that the distance between
$w_4'$ and $y_1$ in $P$ is not more than the distance between
$w_4''$ and $y_1$ in $P$, let $c(w_4w_4')=a$ and
$c(w_4w_4'')=1$ where $w_4', w_4''\in D$.
We will divide three cases to show that $G$ is rainbow connected.

\noindent {\bf Subsubsubcase 2.2.2.1.} $B_1\neq\phi$, $B_2\neq\phi$ and
$B_3\neq\phi$.

First, we can easily check that for each vertex of $A\cup B$, there
is a rainbow path connecting it to any vertex of $P$.

Second, we will show that $G[A_1]$, $G[A_2]$ and $G[A_4]$ are
rainbow connected, respectively. For each pair $(u,v)\in(A_1\times
A_1)$, the edges $uy_1$ and $vv'$ together with the path in $T$
connecting $y_1$ and $v'$ form a rainbow path, where $c(uy_1)=2$
and $c(vv')=a$; For each pair $(u,v)\in(A_2\times A_2)$, the edges
$uy_2$ and $vv'$ together with the path in $T$ connecting $y_2$ and
$v'$ form a rainbow path, where $c(uy_2)=1$ and $c(vv')=a$; For
each pair $(u,v)\in(A_4\times A_4)$, the edges $uu'$ and $vv'$
together with the path in $T$ connecting $u'$ and $v'$ form a
rainbow path, where $c(uu')=1$ and $c(vv')=a$.

Third, we will show that for each vertex of $A_1$, there is a
rainbow path connecting it to any vertex of $A_2\cup A_3\cup A_4\cup
B$. For each pair $(u,v)\in(A_1\times A_2)$, the edges $uu'$ and
$vy_2$ together with the path in $T$ connecting $u'$ and $y_2$ form
a rainbow path, where $c(uu')=a$. For each pair $(u,v)\in(A_1\times
A_3)$, the edges $uu'$ and $vy_1$ together with the path in $T$
connecting $u'$ and $y_1$ form a rainbow path, where $c(uu')=a$.
For each pair $(u,v)\in(A_1\times A_4)$, the edges $uy_1$ and $vv'$
together with the path in $T$ connecting $y_1$ and $v'$ form a
rainbow path, where $c(vv')=a$. For each pair $(u,v)\in(A_1\times
B_1)$, $uy_1v$ is a rainbow path. For each pair $(u,v)\in(A_1\times
B_2)$, $uy_1b_1v$ is a rainbow path, where $b_1\in B_1$. For each
pair $(u,v)\in(A_1\times B_3)$, the edges $uy_1$ and $vv'$ together
with the path in $T$ connecting $y_1$ and $v'$ form a rainbow path.

Fourth, we will show that for each vertex of $A_2$, there is a
rainbow path connecting it to any vertex of $A_3\cup A_4\cup B$. For
each pair $(u,v)\in(A_2\times A_3)$, the edges $uu'$ and $vy_1$
together with the path in $T$ connecting $u'$ and $y_1$ form a
rainbow path, where $c(uu')=a$. For each pair $(u,v)\in(A_2\times
A_4)$, the edges $uy_2$ and $vv'$ together with the path in $T$
connecting $y_2$ and $v'$ form a rainbow path, where $c(vv')=a$.
For each pair $(u,v)\in(A_2\times B_1)$, $uy_2b_2v$ is a rainbow
path, where $b_2\in B_2$. For each pair $(u,v)\in(A_2\times B_2)$,
$uy_2v$ is a rainbow path. For each pair $(u,v)\in(A_2\times B_3)$,
the edges $uy_1$ and $vv'$ together with the path in $T$ connecting
$y_1$ and $v'$ form a rainbow path.

Fifth, we will show that for each vertex of $A_3$, there is a
rainbow path connecting it to any vertex of $A_4\cup B$. For each
pair $(u,v)\in(A_3\times A_4)$, the edges $uy_1$ and $vv'$ together
with the path in $T$ connecting $y_1$ and $v'$ form a rainbow path,
where $c(vv')=a$. For each pair $(u,v)\in(A_3\times B_1)$, $uy_1v$
is a rainbow path. For each pair $(u,v)\in(A_3\times B_2)$, $uy_2v$
is a rainbow path. For each pair $(u,v)\in(A_3\times B_3)$, the
edges $uy_1$ and $vv'$ together with the path in $T$ connecting
$y_1$ and $v'$ form a rainbow path.

Finally, we will show that for each vertex of $A_4$, there is a
rainbow path connecting it to any vertex of $B$. For each pair
$(u,v)\in(A_4\times B_1)$, the edges $uu'$, $vb_2$ and $b_2y_2$
together with the path in $T$ connecting $y_2$ and $u'$ form a
rainbow path, where $c(uu')=1$. For each pair $(u,v)\in(A_4\times
B_2)$, the edges $uu'$ and $vy_2$ together with the path in $T$
connecting $u'$ and $y_2$ form a rainbow path, where $c(uu')=1$.
For each pair $(u,v)\in(A_4\times B_3)$, the edges $uu'$ and $vv'$
together with the path in $T$ connecting $u'$ and $v'$ form a
rainbow path, where $c(uu')=1$.

Hence, we have showed that $G$ is rainbow connected.

From the proof above, we can see the following facts: for any vertex
of $A$ there is a rainbow path connecting it to any vertex of $G$,
and the internal vertex of the rainbow path does not contain any
vertex of $B$; for any vertex of $B_3$, there is a rainbow path
connecting it to any vertex of $G$, and the rainbow path does not
contain any vertex of $B_1\cup B_2$. Hence, in the following proof
we can assume that $B_3=\phi$.

\noindent {\bf Subsubsubcase 2.2.2.2.} $B_1=\phi$ and $B_2\neq\phi$.

We still make use of the above coloring way except for the edges of
$E[A_1, D]$. We now color the edges of $E[A_1, D]$ as follows: For
any vertex $w_1\in A_1$, if $w_1x_1\in E(G)$, then let $c(w_1y_1)=a$
and $c(w_1x_1)=1$; if $w_1x_1\not\in E(G)$, then let $w_1'\in
D\setminus\{y_1, x_1, y_2\}$ with $w_1w_1'\in E(G)$, and let $P[y_1,
w_1']$ be a subpath of $P$, $z\in V(P[y_1, w_1'])$ with
$zw_1'\in E(G)$, then let $c(w_1y_1)=c(zw_1')$,
$c(w_1w_1')=1$. From the coloring, one can easily check that
$G[A_1]$ is rainbow connected.

Now, we show that for each vertex of $A_1$, there is a rainbow path
connecting it to any vertex of $A_2\cup A_3\cup A_4\cup B$. For each
pair $(u,v)\in(A_1\times (A_2\cup A_4))$, the edges $uu'$ and
$vv'$ together with the path in $T$ connecting $u'$ and $v'$ form
a rainbow path, where $c(uu')=1$ and $c(vv')=a$. For each pair
$(u,v)\in(A_1\times A_3)$, $uy_1v$ is a rainbow path. For each pair
$(u,v)\in(A_1\times B_2)$, the edges $uu'$ and $vy_2$ together with
the path in $T$ connecting $u'$ and $y_2$ form a rainbow path,
where $c(uu')=1$.

Thus we have proved that $G$ is rainbow connected.

\noindent {\bf Subsubcase 2.2.2.3} $B_1\neq\phi$ and $B_2=\phi$.

We still make use of the above coloring way except for the edges of
$E[A_2, D]$ and the edges of $E[A_4, D]$. For any vertex $w_4\in
A_4$, we let $c(w_4w_4')=a$ and $c(w_4w_4'')=2$. For any
vertex $w_2\in A_2$, we will color the edges of $E[A_2, D]$ in the
following way: If $w_2x_2\in E(G)$, then let $c(w_2y_2)=a$ and
$c(w_2x_2)=2$; If $w_2x_2\not\in E(G)$, then let $w_2'\in
D\setminus\{y_1, x_2, y_2\}$ with $w_2w_2'\in E(G)$, and let $P[y_2,
w_2']$ be a subpath of $P$, $z'\in V(P[y_2, w_2'])$
with $z'w_2'\in E(G)$, then let $c(w_2y_2)=c(z'w_2')$.
One can easily check that $G[A_2]$ and $G[A_4]$ are rainbow
connected, respectively.

Then, we will show that for each vertex of $A_2$, there is a rainbow
path connecting it to any vertex of $A_1\cup A_3\cup A_4\cup B$. For
each pair $(u,v)\in(A_2\times (A_1\cup A_4))$, the edges $uu'$ and
$vv'$ together with the path in $T$ connecting $u'$ and $v'$ form
a rainbow path, where $c(uu')=2$ and $c(vv')=a$. For each pair
$(u,v)\in(A_2\times A_3)$, $uy_2v$ is a rainbow path. For each pair
$(u,v)\in(A_2\times B_1)$, the edges $uu'$ and $vy_1$ together with
the path in $T$ connecting $u'$ and $y_1$ form a rainbow path,
where $c(uu')=2$.

Finally, we will show that for each vertex of $A_4$, there is a
rainbow path connecting it to any vertex of $A_1\cup A_3\cup B$. For
each pair $(u,v)\in(A_4\times A_1)$, the edges $uu'$ and $vv'$
together with the path in $T$ connecting $u'$ and $v'$ form a
rainbow path, where $c(uu')=2$ and $c(vv')=a$. For each pair
$(u,v)\in(A_4\times A_3)$, the edges $uu'$ and $vy_1$ together with
the path in $T$ connecting $u'$ and $y_1$ form a rainbow path,
where $c(uu')=a$. For each pair $(u,v)\in(A_4\times B_1)$, the
edges $uu'$ and $vy_1$ together with the path in $T$ connecting
$u'$ and $y_1$ form a rainbow path, where $c(uu')=2$.

Thus we have proved that $G$ is rainbow connected.

In the above coloring, we used $e(T)+1$ colors. Hence, $rc(G)\leq
e(T)+1$, and so we can get $rc(G)\leq 2\alpha(G)-1$.

Combining the above Cases 1 and 2, we have completed the proof of Theorem 1. \qed

Since the independence number $\alpha(G)$ is at most the number of
cliques that partition the vertex set of a graph $G$, and the
minimum of such partitions is the chromatic number of the complement
$\bar{G}$ of $G$, we can get the following corollary, which is
Theorem 10 of \cite{Sch1}.
\begin{cor} (Theorem 10, \cite{Sch1} )
Let $G$ be a connected graph with chromatic number $\chi(G)$. Then
$rc(G)\leq 2\chi(\bar{G})-1$.
\end{cor}

\noindent{\bf Acknowledgement.} The authors are very grateful the
referees' for providing us with many new ideas, which helped to
improve our results.

\end{document}